\documentclass[
12pt,
authoryear,
times
]{elsarticle}

\usepackage[colorlinks=true, breaklinks=true 
]{hyperref}
\usepackage{lineno}

\journal{Mathematical Problems in Engineering}

\begin{document}

\begin{frontmatter}

\title{ Efficient formulation of the periodic corrections in Brouwer's gravity solution 
 }

\author
{ Martin Lara\fnref{GRUCACI} 
}
\fntext[GRUCACI]{ Instituto de Ciencia y Technologia, Universidade Federal de S\~ao Paulo }

\ead{mlara0@gmail.com}

\address
{Rua Talim, 330, 12231-280 S\~ao Jos\'e dos Campos, SP, Brazil}


\date{}


\begin{abstract}
The periodic terms of Brouwer's gravity solution are reconstructed in a nonsingular set of variables which are derived from the well-known polar-nodal variables. This change does not affect the essence of the solution, which still keeps all the benefits of the action-angle variables approach, and yields two major improvements. Namely, the periodic corrections of Brouwer's solution are now valid for any eccentricity below one and any inclination except the critical inclination, and, besides, are significantly simpler than the nonsingular corrections in  Lydanne's reformulation of Brouwer's theory. 
\end{abstract}


\end{frontmatter}



\section{Introduction}

Concern in space situational awareness by an increasing number of satellite operators, and in particular the necessity of timely scheduling collision avoidance maneuvers, motivates current interest in improving the capabilities of orbit prediction programs.
\par

Satellite short-term prediction is customarily carried out with SGP4 \citep{HootsRoehrich1980}, an analytical solution which has its roots on Brouwer's celebrated gravity solution to the artificial satellite problem \citep{Brouwer1959}, and that is optimized for the propagation of satellite ephemeris using the element sets in the two-lines format specified by the North American Aerospace Defense Command \citep{HootsSchumacherGlover2004,ValladoCrawfordHujsakKelso2006}. However, it has been claimed that SGP4 may lack of sufficient capabilities for conjunction analysis tasks \citep{KelsoAlfano2005,Kelso2009,HallAlfanoOcampo2010}. Besides, terms that may be missing in SGP4 could be responsible for detected noteworthy along-track errors in the SGP4 predictions for GPS satellites \citep{Kelso2007,Easthope2014}. The known limitations of predicted ephemeris from two-line elements motivate the current development of new algorithms as, for instance, those in the software STELA of the Centre National d'Etudes Spatiales \citep{Fraysseetal2012}.

Brouwer found his solution using a perturbation approach, the so-called von Zeipel method \citep{vonZeipel1916,FerrazMello2007}, which splits the satellite motion into secular terms, long-period corrections, related to the evolution of the argument of the perigee, and short-period corrections, related to the satellite's mean motion. For the secular terms, Brouwer's theory includes gravitational effects up to the second order of $J_2$, the second degree zonal harmonic coefficient of the spherical harmonics expansion of the geopotential, which for the earth is of the order of one thousandth. But in the case of periodic corrections the theory is limited to first order effects of $J_2$. Therefore, the short-period corrections are only related to the contribution of $J_2$, whereas the long-period corrections of Brouwer's gravitational solution include first order corrections due to the few first zonal harmonics \citep{Brouwer1959}.
\par

Brouwer developed his original theory in Delaunay variables, the canonical counterpart of classical Keplerian orbital elements, which, like them, are singular for circular orbits and equatorial orbits. This fact may cause troubles in the computation of the periodic corrections for both low eccentricity and low inclination orbits, but the problem is easily solved by reformulating Brouwer's gravitational solution in nonsingular variables as, for instance, Poincar\'e's canonical variables \citep{Lyddane1963}. However, the periodic corrections either when formulated in Delaunay variables or in Poincar\'e variables are made of long trigonometric series, a fact that very soon motivated efforts in improving their evaluation \citep{Hoots1981}. Among the different efforts in improving the evaluation of Brouwer's gravitational theory, the use of polar-nodal variables for computing the periodic corrections was advocated by different authors \citep{Kozai1962,Izsak1963AJ,Aksnes1972}. These variables, also called Hill or Whittaker variables, are valid for orbits with any eccentricity below 1, but still are singular in the case of equatorial orbits, what may cause trouble in the evaluation of the periodic corrections of almost equatorial orbits.
\par

It worths to remind that Browuer's gravitational theory breaks at the critical inclination of $63.4\deg$. Indeed, since the secular terms in Brouwer's perturbation approach are computed by averaging periodic effects, resonant inclination orbits, and in particular the critical inclination, are excluded from the field of applicability of Brouwer's solution \citep[see][and references therein]{Lara2015IR}.
\par

The efficient evaluation of the periodic corrections is even more critic when taking into account second order corrections, which notably improve the performance of the perturbation theory \citep{Lara2015} but for which the trigonometric series are significantly longer \citep{Kozai1962}, and hence the advantages of using polar-nodal variables are more evident. Besides, the benefits of formulating the periodic corrections in polar-nodal variables are not limited to the case of geopotential perturbations, and this set of canonical variables has revealed very efficient in the evaluation of periodic corrections due to third-body perturbations \citep{LaraVilhenaSanchezPrado2015}.

To avoid the troubles related to the evaluation of the long-period corrections of low-inclination orbits, Aksnes' suggestion of computing the corrections for the satellite's latitude and (true) longitude for these orbits \citep{Aksnes1972} is followed in the present research. Indeed, without limiting to the case of low-inclination orbits, the periodic corrections are rewritten in a set of non-canonical variables which are directly constructed from the polar-nodal ones. These variables are nonsingular and provide a more efficient evaluation of the periodic corrections than the corresponding corrections in Poincar\'e variables which are used in Lyddane's modifications to Brouwer's gravitational solution. 

The paper is organized as follows. For completeness, the construction of short-period corrections in polar-nodal variables of Brouwer's gravitational solution is recalled in Section \ref{s:spe}, while the construction of long-period corrections in polar-nodal variables is illustrated in Section \ref{s:lpe}. Then, the new set of nonsingular, non-canonical variables is introduced in Section \ref{s:cli}, and the long- and short-period corrections are reformulated in the nonsingular variables. The transformations of the nonsingular elements from and to Cartesian variables are free from singularities, and are also documented in Section \ref{s:cli}. 

\section{Short-period elimination} \label{s:spe}

Since this research deals only with perturbations of gravitational origing, the problem of disturbed Keplerian motion can take benefit from Hamiltonian formulation. Thus, the motion of a massless particle in the gravitational field of the earth is derived from the Hamiltonian
\begin{equation}
\mathcal{H}=\mathcal{H}_0+\mathcal{D},
\end{equation}
where $\mathcal{H}_0$ represents the integrable Keplerian Hamiltonian and $\mathcal{D}$ is the disturbing function, which comprises the non-centralities of the geopotential. From the usual solution of Laplace's equation in spherical coordinates, the forces model is further limited to the zonal harmonics case, in which the disturbing function is written
\begin{equation} \label{distpot}
\mathcal{D}=-\frac{\mu}{r}\sum_{m\ge{2}}\left(\frac{\alpha}{r}\right)^mC_{m,0}\,P_{m,0}(\sin\varphi),
\end{equation}
where $\mu$ is the earth's gravitational parameter, $\alpha$ is the earths' equatorial radius, $r$ is the radial distance from the earth's center of mass, $\varphi$ is latitude, $P_{m,0}$ are Legendre polynomials of degree $m$, and $C_{m,0}=-J_m$ are corresponding zonal harmonic coefficients.
\par

The problem of small inclinations in Brouwer's solution is related to the effects of odd zonal harmonics, so to illustrate this case it is enough to consider the impact of $C_{3,0}$, in addition to the main problem, and hence the zonal gravitational potential in Eq.~(\ref{distpot}) is further truncated to the degree $m=3$. Besides, because of the different orders of the harmonic coefficients, where $J_3=\mathcal{O}(J_2)^2$, it is found convenient to make the Hamiltonian perturbative arrangement
\begin{equation} \label{Hj2j3}
\mathcal{H}=H_{0,0}+H_{1,0}+\frac{1}{2}H_{2,0},
\end{equation}
in which
\begin{eqnarray}
H_{0,0} &=& -\frac{\mu}{2a}, \\
H_{1,0} &=& \frac{\mu}{r}\,\frac{1}{4}\,C_{2,0}\,\frac{\alpha^2}{r^2}\left[2-3s^2+3s^2\cos(2f+2\omega)\right], \\
H_{2,0} &=& \frac{\mu}{r}\,\frac{1}{2}\,C_{3,0}\,\frac{\alpha^3}{r^3}\,s\left[
6 \left(1-\frac{5}{4}s^2\right)\sin(f+\omega)
+\frac{5}{2}s^2\sin(3f+3\omega)\right],
\end{eqnarray}
where the relation $\sin\varphi=\sin{I}\sin(f+\omega)$ has been used, with $I$ the orbital inclination, $\omega$ the argument of the perigee, and $f$ the true anomaly, $s$ and $c$ are abbreviations for the sine and cosine of the inclination, respectively, $a$ is the semi-major axis, and, in consequence with the Hamiltonian formulation, all the symbols, that is to say: $a$, $r$, $\omega$, $f$, and $I$, are assumed to be functions of some set of canonical variables.
\par

In particular, Brouwer finds a transformation from ``old'' to ``new'' (or primes) variables, such that the Hamiltonian in the new variables only depends on momenta, whereas the angles have being averaged. Therefore, it relies in the action-angle variables of the Kepler problem, the so-called Delaunay variables. Namely, the mean anomaly $\ell$ and its conjugate momentum $L=\sqrt{\mu\,a}$ (the Delaunay action), the argument of the perigee $g=\omega$ and its conjugate momentum $G=L\,\sqrt{1-e^2}$ (the total angular momentum), where $e$ is the orbital eccentricity, and the argument of the node $h$ and its conjugate momentum $H=G\cos{I}$ (the polar component of the angular momentum). By using this canonical set it is simple to see that $h$ is cyclic in Eq.~(\ref{Hj2j3}) and, therefore, $H$ is an integral of the zonal problem.
\par

The Hamiltonian reduction of Eq.~(\ref{Hj2j3}) by perturbation methods is thoroughly documented in the literature, and hence results are provided without giving details in the method. In particular, the computations carried out were based on the implementation of the Lie transforms method known as Deprit's triangle algorithm, which is nowadays considered standard for Hamiltonian perturbations. Readers interested in the Lie transforms method can find all the required details in the original papers of \citet{Hori1966} and \citet{Deprit1969}, as well as in modern celestial mechanics textbooks like \citep{MeyerHall1992,BoccalettiPucacco1998v2}, or other specialized books as \citep{FerrazMello2007}. Note that, following tradition, in what follows the notation in prime variables is avoided when there is no risk of confusion.
\par

At the first order of the perturbation approach, Deprit's triangle gives
\begin{equation} \label{Dt1st}
H_{0,1}=\left\{H_{0,0},U_1\right\}+H_{1,0},
\end{equation}
where $\{P,Q\}$ notes the Poisson bracket of two functions $P$ and $Q$ of the canonical variables, which in this case are the Delaunay variables. In order to obtain Brouwer's solution, the new Hamiltonian term $H_{0,1}$ is chosen as the averaging of $H_{1,0}$ over the mean anomaly
\begin{equation}
H_{0,1}=-H_{0,0}\,\epsilon_2\,\eta\left(4-6s^2\right),
\end{equation}
where, $\eta$ is the eccentricity function
\begin{equation}
\eta=\frac{G}{L}=\sqrt{1-e^2},
\end{equation}
and, for the sake of abbreviating notation, the function $\epsilon_2\equiv\epsilon_2(G;\mu)$ has been introduced, which is given by
\begin{equation} \label{ep2}
\epsilon_2=\frac{1}{4}C_{2,0}\,\frac{\alpha^2}{p^2},
\end{equation}
where $p=G^2/\mu$ is the \textit{semilatus rectum}.
\par

The corresponding term $U_1$ of the generating function is computed from Eq.~(\ref{Dt1st}) by quadrature
\begin{equation} \label{iU1}
U_1=\frac{1}{n}\int\left(H_{1,0}-H_{0,1}\right)\,\mathrm{d}\ell,
\end{equation}
where $n=\mu^2/L^3$ is the mean motion. In view of the differential relation $a^2\,\eta\,\mathrm{d}\ell=r^2\,\mathrm{d}f$, Eq.~(\ref{iU1}) can be integrated in closed form of the eccentricity to give
\begin{eqnarray} \nonumber
U_1 &=& \frac{1}{2}G\,\epsilon_2\left[(4-6s^2)\,(\phi+e \sin{f}) \right. \\ \label{U1}
&& \left. +3 e s^2 \sin (f+2 g)
+3 s^2 \sin (2 f+2 g)
+e s^2 \sin (3 f+2 g)
\right],
\end{eqnarray}
where $\phi\equiv\phi(\ell,L,G)=f-\ell$ is the equation of the center 
\citep[cf.~Eq.~(15) of][keeping in mind the different sign convention in Hamilton equations]{Brouwer1959}.
\par

Up to the first order, the transformation equations of the averaging are computed from
\begin{equation} \label{teq1}
\rho=\rho'+\{\rho,W_1\},
\end{equation}
where, here, $\rho\in(\ell,g,h,L,G,H)$ and $W_1=U_1$.
\par

Corresponding transformation equations in Delaunay variables can be ex\-press\-ed as Fourier series which involve sine and cosine functions of 10 different arguments of the form $\beta=kf+2mg$ with $k=0,\dots,5$ and $m=-1,0,1$ \citep[cf.~the first order terms in Eqs.~(3.12) and (3.13) of][for instance]{Kozai1962}. 
However, important simplifications can be achieved by using the function
\begin{equation} \label{rRf}
r=\frac{p}{1+e\cos{f}}, 
\end{equation}
instead of wholy expanding the transformation equations as Fourier series. In this way, the number of trigonometric arguments is reduced to just four: $f$, $f+2g$, $2f+2g$, and $3f+2g$ \citep[cf.~Eqs.~(20) and (21) of][]{Brouwer1959}.
\par

Alternatively, as pointed out by \citet{Izsak1963AJ}, the generating function $U_1$ can be expressed in the canonical set of polar-nodal variables $(r,\theta,\nu,R,\Theta,N)$, which stand for the radial distance, the argument of latitude, the argument of the node, the radial velocity, the total angular momentum, and the polar component of the angular momentum. Rewriting Eq.~(\ref{U1}) in polar-nodal variables as $V_1\equiv{U}_1(r,\theta,\nu,R,\Theta,N)$ is straightforward, leading to\footnote{Note that Eq.~(\ref{V1}) differs from Eq.~(5) of \citep{Izsak1963AJ} in the sign. However both equations are equivalent because of the different sign convention used in the derivation of Hamilton equations.}
\begin{equation} \label{V1}
V_1=\epsilon_2\,\frac{\alpha^2}{p^2}\,\Theta
\left[(2-3s^2)\,(\phi+\sigma)+\frac{1}{2} (3+4\kappa)\,s^2 \sin2\theta-\sigma\,s^2\cos2\theta\right],
\end{equation}
where the functions
\begin{equation} \label{ks}
\kappa\equiv\kappa(r,\Theta;\mu)=\frac{p}{r}-1, \qquad \sigma\equiv\sigma(R,\Theta;\mu)=\frac{p\,R}{\Theta},
\end{equation}
are the projections of the eccentricity vector in the orbital frame when written in polar-nodal variables,
which are trivially derived from Eq.~(\ref{rRf}) and its time derivative
\begin{equation}
R=(G/p)\,e\sin{f}.
\end{equation}
\par

The first-order transformation equations of the short-period averaging in polar-nodal variables  are obtained, again, from Eq.~(\ref{teq1}), where, now, $\rho\in(r,\theta,\nu,R,\Theta,N)$ and $W_1=V_1$. In this case, the equation of the center is $\phi\equiv\phi(r,R,\Theta)$, and, in particular, the partial derivatives
\begin{equation}
\frac{\partial\phi}{\partial{r}}=\frac{\sigma}{r}\left(\frac{1+\kappa}{1+\eta}+\frac{\eta}{1+\kappa}\right), \quad
\frac{\partial\phi}{\partial{R}}=\frac{\sigma}{R}\left(\frac{\kappa}{1+\eta}+\frac{2\eta}{1+\kappa}\right), \quad
\frac{\partial\phi}{\partial\Theta}=-\frac{\sigma}{\Theta}\,\frac{2+\kappa}{1+\eta},
\end{equation}
are needed in the computation of the Poisson brackets.
\par

Hence, it is easily obtained
\begin{eqnarray} \label{Dr}
\Delta{r} &=& \epsilon_2\,p \left[\left(2-3 s^2\right) \left(\frac{\kappa}{1+\eta}+\frac{2\eta}{1+\kappa}+1\right)-s^2\cos2\theta\right], 
\\[1ex] \label{Dth}
\Delta\theta &=& \epsilon_2\left\{-3 \left(4-5 s^2\right) \phi
+\left[3-\frac{7}{2}s^2+\left(4-6s^2\right)\kappa\right]\sin2\theta  \right. \\ \nonumber
&& \left.-2 \sigma  \left[5-6s^2+\frac{2+\kappa}{1+\eta}\left(1-\frac{3}{2}s^2\right)+\left(1-2s^2\right)\cos2\theta\right] 
\right\}, 
\\[1ex]  \label{Dnu}
\Delta\nu &=& \epsilon_2\,c \left[6\phi-(4 \kappa +3) \sin2\theta+2\sigma\,(3+\cos2\theta)\right], \\[1ex]
\Delta{R} &=& \epsilon_2\,\frac{\Theta}{p}\left\{2(1+\kappa)^2s^2\sin2\theta-\left(2-3 s^2\right) \sigma 
   \left[\eta+\frac{(1+\kappa)^2}{1+\eta}\right]\right\}, \\[1ex]
\Delta\Theta &=& \epsilon_2\,\Theta\,s^2\left[(3+4\kappa)\cos2\theta+2\sigma\sin2\theta\right], \\[1ex] \label{DN}
\Delta{N} &=& 0,
\end{eqnarray}
which must be evaluated in prime variables for direct corrections $\Delta\rho=\rho-\rho'$ and in original variables for inverse corrections $\Delta\rho'=\rho'-\rho$. Remarkably, now the evaluation of the corrections only requires dealing with sine and cosine functions of the single argument $2\theta$. Note that the evaluation of the equation of the center is required in Eqs.~(\ref{Dth}) and (\ref{Dnu}). It is done using Kepler equation $\phi=f-\ell=f-u+e\sin{u}$, where
\begin{equation}
u=2\arctan\sqrt{\frac{1-e}{1+e}}\tan\frac{f}{2},
\end{equation}
$e=\sqrt{\kappa^2+\sigma^2}$, and $f$ is unambiguously computed from $\cos{f}=\kappa/e$ and $\sin{f}=\sigma/e$.
\par

The second order of Deprit's triangle gives 
\begin{equation} \label{Dt2nd}
H_{0,2}=\left\{H_{0,0},U_2\right\}+\left\{H_{0,1},U_1\right\}+\left\{H_{1,0},U_1\right\}+H_{2,0},
\end{equation}
and the new Hamiltonian term $H_{0,2}$ is chosen as the average of $H_{2,0}$ plus the computable Poisson brackets, to give
\begin{eqnarray}
H_{0,2} &=& H_{0,0}\left\{\frac{3}{2}\epsilon_2^2
\left[5(8-16s^2+7s^4)+\eta\,(4-6s^2)^2-\eta^2\,(8-8s^2-5s^4) \right. \right. \qquad \\ \nonumber
&& \left.\left. -2(14-15s^2)\,s^2\,e^2\cos2g\right]\eta
-\frac{3}{2}C_{3,0}\,\frac{\alpha^3}{p^3}\,(4-5s^2)s\,\eta\,e\sin{g} \right\}.
\end{eqnarray}
Second order corrections to the orbital elements are of the order of the square of $J_2$, and are normally omitted. Therefore, there is no need of solving $U_2$ from Eq.~(\ref{Dt2nd}), and the short-period transformation limits to the first order corrections in Eqs.~(\ref{Dr})--(\ref{DN}), whose simple inspection shows that are free from singularities either for equatorial or circular orbits.
\par


\section{Long-period elimination} \label{s:lpe}

After the mean anomaly has been averaged, the long-period Hamiltonian is
\begin{equation} \label{Hpj2j3}
\mathcal{K}=K_{0,0}+K_{1,0}+\frac{1}{2}K_{2,0}
\end{equation}
where $K_{0,0}=H_{0,0}$, $K_{1,0}=H_{0,1}$, $K_{2,0}=H_{0,2}$, which are expressed in prime elements, although primes have been dropped for alleviating notation.
\par

In the new notation, the first order of Deprit's triangle in Eq.~(\ref{Dt1st}) is rewritten
\begin{equation} \label{Dt1stp}
K_{0,1}=\left\{K_{0,0},X_1\right\}+K_{1,0}.
\end{equation}
Because $K_{1,0}$ in Eq.~(\ref{Hpj2j3}) does not depend on $g$, the new first-order Hamiltonian term is chosen $K_{0,1}=K_{1,0}$, and hence $\{K_{0,0},X_1\}=0$ from Eq.~(\ref{Dt1stp}). However, this does not mean to make null the first order term $X_1$ of the long-period generating function. Quite on the contrary, since the generating function of the long-period averaging does not depend on $\ell$ then $\{K_{0,0},X_1\}$ necessarily vanishes in Eq.~(\ref{Dt1stp}). Therefore, the term $X_1$ can only be determined at the next order of the perturbation algorithm.
\par

The Poisson bracket $\{K_{0,0},X_2\}$ vanishes likewise, and the second order of Deprit's triangle in Eq.~(\ref{Dt2nd}) is simplified in this case to
\begin{equation} \label{Dt2ndp}
K_{0,2}=2\left\{K_{0,1},X_1\right\}+K_{2,0},
\end{equation}
where the term $K_{0,2}$ is chosen as the average of $K_{2,0}$ over the argument of the perigee. That is,
\begin{equation}
K_{0,2}=K_{0,0}\frac{3}{2}\epsilon_2^2\,\eta\left[
5\left(8-16s^2+7s^4\right)
+\left(4-6s^2\right)^2\eta
-\left(8-8s^2-5s^4\right)\eta^2\right].
\end{equation}
\par

It follows the computation of $X_1$ from Eq.~(\ref{Dt2ndp}) by a quadrature:
\begin{equation}
X_1=\frac{1}{n}\,\frac{2 a^2 \eta ^4}{3 \alpha ^2\left(4-5 s^2\right) C_{2,0}}\,\int \left(H_{0,2}-H_{2,0}\right) \,\mathrm{d}g,
\end{equation}
which is trivially solved to give
\begin{equation} \label{X1}
X_1 = G \left(-\epsilon_2\,\frac{14-15s^2}{4-5s^2}\,\frac{1}{8}s^2\,e^2\sin2g+\epsilon_3\,s\,e\cos{g}\right), 
\end{equation}
where, for conciseness, the notation
\begin{equation} \label{ep3}
\epsilon_3 = \frac{1}{2}\frac{\alpha}{p}\,\frac{C_{3,0}}{C_{2,0}},
\end{equation}
has been introduced.
\par

Next, the long-period generating function
\begin{eqnarray} \label{Y1}
Y_1 &=& -\epsilon_2\,\Theta\,s^2\frac{14-15s^2}{8\left(4-5s^2\right)}\left[
\left(\kappa^2-\sigma^2\right)\sin2\theta-2\kappa\,\sigma\cos2\theta\right] \\[1ex]  \nonumber
&& + \epsilon_3\,\Theta\,s\,(\kappa\cos\theta+\sigma\sin\theta).
\end{eqnarray}
is obtained by rewriting Eq.~(\ref{X1}) in polar-nodal variables, and the first order transformation equations in polar-nodal variables are obtained again from Eq.~(\ref{teq1}), whit $\rho\in(r,\theta,\nu,R,\Theta,N)$ and where now $W_1$ is replaced by $Y_1$. In this way the long-period corrections
\begin{eqnarray} \label{dr}
\delta{r} &=& 
p\left[\epsilon_2\,s^2\,\frac{1-15c^2}{4\left(1-5c^2\right)}\,(\kappa\cos2\theta+\sigma\sin2\theta)
+\epsilon_3\,s\sin\theta\right], 
\\[1ex] \label{dtheta}
\delta\theta &=& 
\epsilon_2\frac{1}{2\left(1-5c^2\right)^2}\left[(q_2+q_5\,\kappa)\sigma\cos2\theta
  -(q_1\,\sigma^2 +q_2\,\kappa+q_3\,\kappa^2)\sin2\theta
   \right] \\ \nonumber
&& +\epsilon_3\left[\left(\frac{\kappa }{s}+2s\right)\cos\theta+\left(\frac{1}{s}-s\right)\sigma\sin\theta\right], 
\\[1ex] \label{dnu}
\delta\nu &=&
\epsilon_2\,\frac{q_6}{4\left(1-5c^2\right)^2}
\left[\left(\kappa^2-\sigma^2\right)\sin2\theta-2\kappa\,\sigma\cos2\theta\right] \\ \nonumber
&& -\epsilon_3\,\frac{c}{s}\left(\kappa\cos\theta+\sigma\sin\theta\right), \\[1ex]
\delta{R} &=&
\frac{\Theta}{p}\,(1+\kappa)^2\left[
\epsilon_2\,\frac{1-15c^2}{4\left(1-5c^2\right)}\,s^2\,(\sigma\cos2\theta-\kappa\sin2\theta)
+\epsilon_3\,s\cos\theta
\right], \\[1ex] \label{deZ}
\delta\Theta &=&
\Theta\,\epsilon_2\,\frac{1-15c^2}{4\left(1-5c^2\right)}\,s^2\left[(\kappa^2-\sigma^2)\cos2\theta+2\kappa\,\sigma\sin2\theta\right] \\ \nonumber
&& +\Theta\,\epsilon_3\,s\,(\kappa\sin\theta-\sigma\cos\theta), \\[1ex] \label{deN}
\delta{N} &=& 0,
\end{eqnarray}
have been obtained, where the inclination polynomials $q_j$, $j=1,\dots6$ are given in Table \ref{t:iq}. Equations (\ref{dr})--(\ref{deN}) must be evaluated in second prime variables for direct corrections $\Delta\rho'=\rho'-\rho''$ and in prime variables from the inverse corrections $\Delta\rho''=\rho''-\rho'$.
\par

\begin{table*}[htb]
\caption{Inclination polynomials.} 
\centering
\begin{tabular}{@{}rl@{}}\hline
$q_{0}=$ & $(1-15c^2)\,(1-5c^2)$ \vphantom{$\left(\frac{M}{M}\right)^M$} \\[0.5ex]
$q_{1}=$ & $\frac{1}{4}(1-43c^2+155c^4-225c^6)$ \\[0.5ex]
$q_{2}=$ & $s^2\,q_0$  \\[0.5ex]
$q_{3}=$ & $\frac{1}{4}(1+c^2+35c^4+75c^6)$ \\[0.5ex]
$q_{5}=$ & $c^2\,(11-30c^2+75c^4)$ \\[0.5ex]
$q_{6}=$ & $q_5/c$ \\[0.5ex]
$q_{7}=$ & $\frac{1}{4}(1+3c^2-5c^4+225c^6)$ \\[0.5ex]
$q_{8}=$ & $\frac{1}{4}(1-45c^2+195c^4-375c^6)$ \\[0.5ex]
$q_{9}=$ & $\frac{1}{4}(1+75c^4)$ \\[0.5ex]
$q_{10}=$ & $\frac{1}{4}(1-40c^2+75c^4)$ \\[0.5ex]
$q_{11}=$ & $2c^2\,(6-25c^2+75c^4)$ \\[0.5ex]
$q_{12}=$ & $10c^2$ \\[0.5ex]
$q_{13}=$ & $q_{0}\,(1+c)$ \\[0.5ex]
$q_{14}=$ & $\frac{1}{4}(1-c)\,(1-20c-40c^2+75c^4)$ \\[0.5ex]
$q_{15}=$ & $\frac{1}{4}(1+23c-20c^2-80c^3+75c^4+225c^5)$ \\[0.5ex]
\hline
\end{tabular}
\label{t:iq}
\end{table*}

Note that the term $1-5c^2$ in denominators of Eqs.~(\ref{dr})--(\ref{deZ}) prevents application of the long-period corrections to orbits with the critical inclination of $63.4\deg$. This singularity is not related to the variables used, and simply reflects the fact that inclination resonances are out of the range of applicability of Brouwer's gravitational solution \citep[see][and references therein]{Lara2015IR}.

Finally, it deserves mentioning that after computing the double-prime Delaunay variables from the secular terms, the Kepler equation must be solved to find first $f$, and then $\theta$, in order to compute corresponding double-prime polar-nodal variables.

\section{The case of low inclinations} \label{s:cli}

Due to the contribution of the odd zonal harmonic $C_{3,0}$, it happens that $\delta\theta$ and $\delta\nu$ are singular for equatorial orbits. However, as the simple inspection of Eqs.~(\ref{dtheta}) and (\ref{dnu}) may suggest, the trouble in the case of low inclinations is easily remedied by computing the long-period corrections to the nonsingular, non-canonical elements $(\psi,\xi,\chi,r,R,\Theta)$, where
\begin{equation} \label{ycx}
\psi=\theta+\nu, \qquad \xi=s\sin\theta, \qquad \chi=s\cos\theta.
\end{equation}
\par
Indeed, by simple differentiation
\begin{eqnarray} \label{deltapsi}
\delta\psi &=& \delta\theta+\delta\nu, \\ \label{deltaxi}
\delta\xi &=& \left(\frac{\delta\Theta}{s}\right)\frac{c^2}{\Theta}\sin\theta+(s\,\delta\theta)  \cos\theta, \\ \label{deltachi}
\delta\chi &=& \left(\frac{\delta\Theta}{s}\right)\frac{c^2}{\Theta}\cos\theta-(s\,\delta\theta)\sin\theta,
\end{eqnarray}
where neither $(\delta\Theta/s)$ nor $(s\,\delta\theta)$ are affected of singularities, as easily checked in Eqs.~(\ref{deZ}) and (\ref{dtheta}), respectively. Alternatively, because Eq.~(\ref{teq1}) applies to any function of the canonical variables \citep{Deprit1969,DepritRom1970}, the corrections in Eqs.~(\ref{deltapsi})--(\ref{deltachi}), can be computed from corresponding Poisson brackets, viz.
\begin{eqnarray}
\delta\psi &=& \{\theta+\nu,Y_1\} \\ 
\delta\xi &=& \{s\sin\theta,Y_1\}, \\ 
\delta\chi &=& \{s\cos\theta,Y_1\}.
\end{eqnarray}
Straightforward manipulations lead to the explicit equations in the nonsingular variables
\begin{eqnarray} \label{dyns}
\delta\psi \!\!\!&=&\!\!\! \frac{1}{1+c}\,\Bigg\{\epsilon_2\frac{1}{2(1-5c^2)}\left[
2\xi\,\chi\,(q_{13}\,\kappa + q_{14}\,\kappa^2 + q_{15}\,\sigma^2) \right. \\ \nonumber
&& \left.-\sigma\,(\chi^2-\xi^2)\,(q_{13} - q_6\,\kappa)\right]
+\epsilon_3\,\Big[(2+2c+\kappa)\,\chi-c\,\sigma\,\xi\Big] \Bigg\},
\\[1ex]  \label{dxins}
\delta\xi \!\!\!&=&\!\!\! \frac{\epsilon_2}{4(1-5c^2)}\left[P_1\,\xi+P_2\,(3\chi^2-\xi^2)\,\xi-P_3\,\sigma\,\chi-P_4\,\sigma\,(\chi^2-3\xi^2)\,\xi\right] \\ \nonumber
&& +\frac{1}{2}\epsilon_3\left[2s^2+(1+c^2)\,\kappa+(2+\kappa)\,(\chi^2-\xi^2)\right],
\\[1ex] \label{dchins}
\delta\chi \!\!\!&=&\!\!\! \!\!-\frac{\epsilon_2}{4(1-5c^2)}\left[P_1\,\chi+P_2\,(3\xi^2-\chi^2)\,\chi+P_3\,\sigma\,\xi+P_4\,\sigma\,(\xi^2-3\chi^2)\,\xi\right] \\ \nonumber
&& -\epsilon_3\left[c^2\,\sigma +(2+\kappa)\,\chi\,\xi\right],
\\[1ex] \label{drns}
\delta{r} \!\!\!&=&\!\!\! \epsilon_2\frac{1-15c^2}{4(1-5c^2)}
\left[2\sigma\,\xi\,\chi -\kappa\,(\xi^2-\chi^2)\right]
+\epsilon_3\,\xi,
\\[1ex] \label{dRRns}
\delta{R} \!\!\!&=&\!\!\! \frac{\Theta}{p}(1+\kappa)^2 \left\{
-\epsilon_2\frac{1-15c^2}{4(1-5c^2)}\left[
2\kappa\,\xi\,\chi+\sigma\,(\xi^2-\chi^2)\right]
+\epsilon_3\,\chi
\right\},
\\[1ex] \label{dZZns}
\delta\Theta \!\!\!&=&\!\!\!
\Theta\left\{\epsilon_2\,\frac{1-15c^2}{4\left(1-5c^2\right)}\left[(\kappa^2-\sigma^2)\,(\chi^2-\xi^2)+4\kappa\,\sigma\,\chi\,\xi\right] + \epsilon_3(\kappa\,\xi-\sigma\,\chi)\right\}, 
\end{eqnarray}
where the coefficients $P_j$ ($j=1,\dots4$) are given in Table \ref{t:iep}, $s^2=\xi^2+\chi^2$ from Eq.~(\ref{ycx}), $c=\sqrt{1-s^2}$, and $\kappa$ and $\sigma$ are given in Eq.~(\ref{ks}). Note the almost symmetric form of the corrections $\delta\xi$ and $\delta\chi$ in Eqs.~(\ref{dxins}) and (\ref{dchins}), respectively.

%
\begin{table*}[htbp]
\caption{Coefficients $P_j$ in Eqs. (\protect\ref{dxins})--(\protect\ref{dchins}).}
\centering
\begin{tabular}{@{}rl@{}}\hline
$P_{1}=$ & $q_2\,\kappa + q_7\,\kappa^2 + q_8\,\sigma^2$ \vphantom{$\left(\frac{M}{M}\right)^M$} 
\\[0.5ex]
$P_{2}=$ & $q_{0}\,\kappa+q_9\,\kappa^2 + q_{10}\,\sigma^2$ \\[0.5ex]
$P_{3}=$ & $q_2+ q_{11}\,\kappa$  \\[0.5ex]
$P_{4}=$ & $q_{0} + q_{12}\,\kappa$ \\[0.5ex]
\hline
\end{tabular}
\label{t:iep}
\end{table*}
%

In the case of the earth, $J_2\approx\sin^22^\circ$, and hence terms of the order of $s^2$ and higher can be neglected for the lower inclination orbits, because they only produce higher order effects. Therefore, the corrections
\begin{eqnarray} \label{dy}
\delta\psi &=& \epsilon_3\,\frac{1}{2}\left[\chi\,(4+\kappa) -\xi\,\sigma\right], \\[1ex] \label{dxi}
\delta\xi &=& \epsilon_2\,\frac{7}{8}\left[(\kappa^2-\sigma^2)\,\chi+2\kappa\,\sigma\,\xi\right]
-\epsilon_3\,\sigma, \\[1ex] \label{dchi}
\delta\chi &=& -\epsilon_2\,\frac{7}{8}\left[(\kappa^2-\sigma^2)\,\xi-2\kappa\,\sigma\,\chi\right]
+\epsilon_3\,\kappa, \\[1ex]
\delta{r} &=& 
\epsilon_3\,\xi\,p, \\[1ex]
\delta{R} &=&
\epsilon_3\,(1+\kappa)^2\,\chi\,\frac{\Theta}{p}, \\[1ex]  \label{dZ}
\delta\Theta &=&
\epsilon_3\,(\kappa\,\xi-\sigma\,\chi)\,\Theta,
\end{eqnarray}
are straightforwardly derived from Eqs.~(\ref{dr})--(\ref{deN}). Note that Eqs.~(\ref{dy})--(\ref{dchi}) have been previously provided by \citet{Aksnes1972}.
\par

\subsection{Transformation from Cartesian variables}
The direct transformation from nonsingular to Cartesian variables is obtained by means of the usual rotations applied to the projections of the position and velocity vectors in the orbital frame. Thus,
\begin{equation} \label{rotations}
\left(\begin{array}{cc} x & X \\ y & Y \\ z & Z \end{array}\right)=
R_3(-\nu)\circ{R}_1(-I)\circ{R}_3(-\theta)\circ\left(\begin{array}{cc} r & \phantom{\theta}\dot{r}=R\quad \\ 0 & r\,\dot\theta=\Theta/r \\ 0 & 0 \end{array}\right),
\end{equation}
where $R_1$, $R_3$, are the usual rotation matrices
\begin{equation}
R_1(\beta)=\left(\begin{array}{ccc} 1 & 0 & 0 \\ 0 & \cos\beta & \sin\beta \\ 0 & -\sin\beta & \cos\beta \end{array}\right), \quad
R_3(\beta)=\left(\begin{array}{ccc} \cos\beta & \sin\beta & 0 \\ -\sin\beta & \cos\beta & 0 \\ 0 & 0 & 1\end{array}\right).
\end{equation}

\par

After replacing $\nu=\psi-\theta$ and $\sin\theta=\xi/s$, $\cos\theta=\chi/s$, in Eq.~(\ref{rotations}), the transformation from nonsingular to Cartesian variables can be obtained from the sequence
\begin{eqnarray} \label{ns2x}
x &=& r\,(t\cos\psi+q\sin\psi), \\ \label{ns2y}
y &=& r\,(t\sin\psi-q\cos\psi), \\
z &=& r\,\xi, \\
X &=& R\,(t\cos\psi+q\sin\psi) - \frac{\Theta}{r}\,(q\cos\psi+\tau\sin\psi), \\ \label{ns2Y}
Y &=& R\,(t\sin\psi-q\cos\psi) - \frac{\Theta}{r}\,(q\sin\psi-\tau\cos\psi), \\  \label{ns2zz}
Z &=& R\,\xi+\frac{\Theta}{r}\,\chi,
\end{eqnarray}
where 
\begin{equation} \label{tq}
t=1 - \frac{\xi^2}{1 + c}, \qquad
\tau=1 - \frac{\chi^2}{1 + c}, \qquad
q=\frac{\xi\,\chi}{1 + c},
\end{equation}
and $c=N/\Theta$. Remark that $N$ is an integral of the zonal problem and, therefore, his value is always known from given initial conditions.
\par

The inverse transformation, from Cartesian to nonsingular variables, is obtained from the sequence
\begin{eqnarray} \label{rx}
r &=& \sqrt{x^2+y^2+z^2}, \\[1ex] \label{Rx}
R &=& \frac{1}{r}\left(x\,X+y\,Y+z\,Z\right), \\[1ex] \label{Nx}
N &=& x\,Y-y\,X, \\[1ex] \label{Zx}
\Theta &=& \sqrt{(y\,Z-z\,Y)^2+(z\,X-x\,Z)^2+N^2}, \\[1ex]
\chi &=& \frac{1}{\Theta}\left(r\,Z-z\,R\right), \\[1ex]
\xi &=& \frac{z}{r}, \\[1ex] \label{scy}
\sin\psi &=& \frac{x\,q+y\,t}{(t^2+q^2)\,r}, \qquad \cos\psi \;=\; \frac{x\,t-y\,q}{(t^2+q^2)\,r},
\end{eqnarray}
where the computation of $\psi$, which is unambiguously determined from Eq.~(\ref{scy}), requires the previous computation of $t$ and $q$ from Eq.~(\ref{tq}).
\par

Note that Eqs.~(\ref{tq}) are singular for equatorial retrograde orbits, a case in which $c=-1$. However, this drawback is easily remedied and the case of almost equatorial, retrograde orbits is effectively addressed by using the variable $\psi^*=\theta-\nu$ instead of $\psi$. Then, the corrections in Eqs.~(\ref{dy})--(\ref{dZ}) still apply, yet the conversion from nonsingular to Cartesian coordinates is slightly modified. Indeed, $y$ and $Y$ in Eqs.~(\ref{ns2y}) and (\ref{ns2Y}) must be replaced by $-y$ and $-Y$, respectively, whereas changing $c$ in Eq.~(\ref{tq}) by $|c|$ allows for computing $t$ and $q$ from this equation in both cases of direct and retrograde inclinations.
\par

\subsection{Short-period corrections in nonsingular variables}

In spite of there is no trouble in the evaluation of the short-period corrections in the case of low-inclination orbits, it may be convenient to compute Eqs.~(\ref{Dr})--(\ref{DN}) also in nonsingular variables. In this case, 
\begin{eqnarray}
\Delta\psi &=& \epsilon_2\left\{(3+6c-15c^2)\,\phi
+\sigma\left[2+6 c-12c^2+(1-3c^2)\frac{2+\kappa}{1+\eta}
\right.\right. \\ \nonumber
&& \left.\left. +\frac{2+4c}{1+c}\,(\chi ^2-\xi ^2)\right]
-\frac{1+7c+4(1+3c)\kappa}{1+c}\,\xi\,\chi \right\},
\\
\Delta\xi &=& \epsilon_2 \left\{
\sigma\left[4\chi^2-12c^2+(1-3c^2)\frac{2+\kappa}{1+\eta}\right]\chi \right. \\ \nonumber
&& \left. -\left[(1+4\kappa)\,\chi^2 - (3+4\kappa)\,c^2\right]\xi
+3(1-5c^2)\,\phi\,\chi
\right\},
\\[1ex]
\Delta\chi &=& -\epsilon_2 \left\{ 
\sigma\left[4\chi^2-8c^2+(1-3c^2)\frac{2+\kappa}{1+\eta}\right]\xi \right. \\ \nonumber
&& \left. -\left[(1+4\kappa)\,\xi^2 - (3+4\kappa)\,c^2\right]\chi
+3(1-5c^2)\,\phi\,\xi
\right\},
\\
\Delta{r} &=& \epsilon_2\,p
\left[\xi ^2-\chi ^2+\left(1+\frac{\kappa }{1+\eta}+\frac{2 \eta }{1+\kappa}\right)(2-3s^2)\right],
\\
\Delta{R} &=& \epsilon_2
\frac{\Theta}{p}\left[4(1+\kappa)^2\,\xi\,\chi
-\sigma\left(\eta+\frac{(1+\kappa)^2}{1+\eta}\right)
(2-3s^2)\right],
\\
\Delta\Theta &=&  \epsilon_2\,\Theta\left[(3+4\kappa)\,(\xi^2-\chi^2)-4\sigma\,\xi\,\chi\right],
\end{eqnarray}
and, except for higher order effects, the short-period corrections to the lower inclination orbits may be written in nonsingular elements as
\begin{eqnarray} \label{Dpsi0}
\Delta\psi &=& -2\epsilon_2\left[3\phi+\left(2+\frac{2+\kappa}{1+\eta}\right)\sigma\right], \\[1ex] \label{Dxi0}
\Delta\xi  &=&  \epsilon_2\left[(3+4\kappa)\,\xi-2\left(6+\frac{2+\kappa}{1+\eta}\right)\sigma\,\chi-12\phi\,\chi\right], \\ \label{Dchi0}
\Delta\chi &=& -\epsilon_2\left[(3+4\kappa)\,\chi-2\left(4+\frac{2+\kappa}{1+\eta}\right)\sigma\,\xi-12\phi\,\xi\right], \\
\Delta{r} &=& 2\epsilon_2\,p \left(1+\frac{\kappa}{1+\eta}+\frac{2\eta}{1+\kappa}\right), \\[1ex]
\Delta{R} &=& -2\epsilon_2\,\frac{\Theta}{p}\,\sigma\left[\eta+\frac{(1+\kappa)^2}{1+\eta}\right], \\[1ex] \label{DeltaZ0}
\Delta\Theta &=& 0.
\end{eqnarray}

Finally, in spite of a state $(x,y,0,X,Y,0)$ corresponding to an exactly (instantaneous) equatorial orbit would be rarely obtained when working in real arithmetic, it deserves to mention that in the space (original or double prime) that it might happen 
$\xi=\chi=0$ and it does not make sense to speak of the node or the argument of latitude. However, periodic corrections still exist for  $\xi$ and $\chi$. Indeed, while short-period corrections $\Delta\xi$ and $\Delta\chi$ vanish for equatorial orbits, as derived from Eqs~(\ref{Dxi0})--(\ref{Dchi0}), 
corresponding long-period corrections do not, and Eqs.~(\ref{dxi}) and (\ref{dchi}) result in
\begin{equation}
\delta\xi= -\epsilon_3\,\sigma, \qquad
\delta\chi=\epsilon_3\,\kappa.
\end{equation}

\section{Conclusions}

Soon after Brouwer's solution was announced, the reformulation in polar-nodal variables of both the short- and long-period corrections was suggested as a way of simplifying their evaluation. Indeed, as odd as it may seem to introduce short-period terms in the computation of long-period corrections, this artifact prevents the usual deterioration of the corrections in the case of low-eccentricity orbits, yet the case of low-inclination orbits must be treated separately. However, the elementary inspection of the long-period corrections in polar-nodal variables reveals a simple set of (non-canonical) elements that may be used for dealing properly with that case. The new formalism is nonsingular, yields significantly less computational effort than Lyddane's nonsingular variables approach, and can be extended to reformulate third-body periodic corrections in a compact form. The latter is under development and will be published elsewhere. 


%
%

\section*{Acknowledgemnts}

Partial support from project ESP2013-41634-P of the Ministry of Economic Affairs and Competitiveness of Spain is recognized.


\section*{References}
\begin{thebibliography}{25}
\expandafter\ifx\csname natexlab\endcsname\relax\def\natexlab#1{#1}\fi
\expandafter\ifx\csname url\endcsname\relax
  \def\url#1{\texttt{#1}}\fi
\expandafter\ifx\csname urlprefix\endcsname\relax\def\urlprefix{URL }\fi

\bibitem[{{Aksnes}(1972)}]{Aksnes1972}
{Aksnes}, K., Feb. 1972. {On the Use of the Hill Variables in Artificial
  Satellite Theory}. Astronomy and Astrophysics 17, 70--75.

\bibitem[{{Boccaletti} and {Pucacco}(2002)}]{BoccalettiPucacco1998v2}
{Boccaletti}, D., {Pucacco}, G., 2002. {Theory of orbits. Volume 2:
  Perturbative and geometrical methods}, 1st Edition. Astronomy and
  Astrophysics Library. Springer-Verlag, Berlin Heidelberg New York.

\bibitem[{Brouwer(1959)}]{Brouwer1959}
Brouwer, D., November 1959. Solution of the problem of artificial satellite
  theory without drag. The Astronomical Journal 64, 378--397.

\bibitem[{{Deprit}(1969)}]{Deprit1969}
{Deprit}, A., 1969. Canonical transformations depending on a small parameter.
  Celestial Mechanics 1~(1), 12--30.

\bibitem[{Deprit and Rom(1970)}]{DepritRom1970}
Deprit, A., Rom, A., June 1970. {The Main Problem of Artificial Satellite
  Theory for Small and Moderate Eccentricities}. Celestial Mechanics 2~(2),
  166--206.

\bibitem[{{Easthope}(2014)}]{Easthope2014}
{Easthope}, P.~F., 2014. {Examination of SGP4 along-track errors for initially
  circular orbits}. IMA Journal of Applied Mathematics on line, 1--15.

\bibitem[{{Ferraz-Mello}(2007)}]{FerrazMello2007}
{Ferraz-Mello}, S., Jan. 2007. {Canonical Perturbation Theories - Degenerate
  Systems and Resonance}. Vol. 345 of Astrophysics and Space Science Library.
  Springer, New York.

\bibitem[{{Fraysse} et~al.(2012){Fraysse}, {Le Fevre}, {Morand}, {Deleflie},
  {Mercier}, and {Dental}}]{Fraysseetal2012}
{Fraysse}, H., {Le Fevre}, C., {Morand}, V., {Deleflie}, F., {Mercier}, P.,
  {Dental}, C., 2012. {STELA, a tool for long term orbit propagation}. In:
  Proceedings of the 5th International Conference on Astrodynamics Tools and
  Techniques. ESTEC/ESA, The Netherlands.

\bibitem[{Hall et~al.(2010)Hall, Alfano, and Ocampo}]{HallAlfanoOcampo2010}
Hall, R., Alfano, S., Ocampo, A., 2010. {Advances in Satellite Conjunction
  Analysis}. In: Presented at tne 2010 AMOS Conference. USA.

\bibitem[{{Hoots}(1981)}]{Hoots1981}
{Hoots}, F.~R., Aug. 1981. {Reformulation of the Brouwer geopotential theory
  for improved computational efficiency}. Celestial Mechanics 24, 367--375.

\bibitem[{{Hoots} and {Roehrich}(1980)}]{HootsRoehrich1980}
{Hoots}, F.~R., {Roehrich}, R.~L., December 1980. {Models for Propagation of
  the NORAD Element Sets}. {Project SPACETRACK, Rept.~3, U.S.~Air Force
  Aerospace Defense Command, Colorado Springs, CO}.

\bibitem[{{Hoots} et~al.(2004){Hoots}, {Schumacher}, and
  {Glover}}]{HootsSchumacherGlover2004}
{Hoots}, F.~R., {Schumacher}, Jr., P.~W., {Glover}, R.~A., Mar. 2004. {History
  of Analytical Orbit Modeling in the U. S. Space Surveillance System}. Journal
  of Guidance, Control, and Dynamics 27~(5), 174--185.

\bibitem[{{Hori}(1966)}]{Hori1966}
{Hori}, G., 1966. {Theory of General Perturbation with Unspecified Canonical
  Variables}. Publications of the Astronomical Society of Japan 18~(4),
  287--296.

\bibitem[{{Izsak}(1963)}]{Izsak1963AJ}
{Izsak}, I.~G., Oct. 1963. {A note on perturbation theory}. The Astronomical
  Journal 68, 559--561.

\bibitem[{Kelso(2007)}]{Kelso2007}
Kelso, T.~J., 2007. {Validation of SGP4 and IS-GPS-200D against GPS precision
  ephemerides (AAS 07-127)}. In: Proceedings of the 17th AAS/AIAA Spaceflight
  Mechanics Conference, Sedona, AZ. American Astronautical Society, Univelt,
  Inc., USA.
\newline\urlprefix\url{https://celestrak.com/publications/AAS/07-127/AAS-07-127.pdf}

\bibitem[{Kelso(2009)}]{Kelso2009}
Kelso, T.~J., 2009. {Analysis of the Iridium 33--Cosmos 2251 Collision (AAS
  09-368)}. In: Paper AAS 09-368, American Astronautical Society. American
  Astronautical Society, Univelt, Inc., USA.

\bibitem[{Kelso and Alfano(2005)}]{KelsoAlfano2005}
Kelso, T.~J., Alfano, S., 2005. {Satellite Orbital Conjunction Reports
  Assessing Threatening Encounters in Space (SOCRATES) (AAS 05-124)}. In:
  Proceedings of the AAS/AIAA Spaceflight Mechanics Conference, Copper
  Mountain, CO. American Astronautical Society, Univelt, Inc., USA.

\bibitem[{{Kozai}(1962)}]{Kozai1962}
{Kozai}, Y., September 1962. {Second-Order Solution of Artificial Satellite
  Theory without Air Drag}. The Astronomical Journal 67~(7), 446--461.

\bibitem[{Lara(2015)}]{Lara2015}
Lara, M., 2015. {LEO intermediary propagation as a feasible alternative to
  Brouwer's gravity solution}. Advances in Space Research on line~(December
  2014).
\newline\urlprefix\url{http://dx.doi.org/10.1016/j.asr.2014.12.023}

\bibitem[{{Lara}({2015b})}]{Lara2015IR}
{Lara}, M., {2015b}. {On Inclination Resonances in Artificial Satellite
  Theory}. Acta Astronautica on line.
\newline\urlprefix\url{http://dx.doi.org/10.1016/j.actaastro.2015.02.0013}

\bibitem[{Lara et~al.(2015)Lara, Vilhena~de Moraes, Sanchez, and
  Prado}]{LaraVilhenaSanchezPrado2015}
Lara, M., Vilhena~de Moraes, R., Sanchez, D.~M., Prado, A. F. B.~A., 2015.
  {Efficient computation of short-period analytical corrections due to
  third-body effects (AAS 15-295)}. American Astronautical Society, Univelt,
  Inc., USA.

\bibitem[{{Lyddane}(1963)}]{Lyddane1963}
{Lyddane}, R.~H., Oct. 1963. {Small eccentricities or inclinations in the
  Brouwer theory of the artificial satellite}. {Astronomical Journal} 68,
  555--558.

\bibitem[{{Meyer} and {Hall}(1992)}]{MeyerHall1992}
{Meyer}, K.~R., {Hall}, G.~R., 1992. {Introduction to Hamiltonian Dynamical
  Systems and the N-Body Problem}. Springer, New York.

\bibitem[{Vallado et~al.(2006)Vallado, Crawford, Hujsak, and
  Kelso}]{ValladoCrawfordHujsakKelso2006}
Vallado, D.~A., Crawford, P., Hujsak, R., Kelso, T.~S., August 2006.
  {Revisiting Spacetrack Report \#3 (AIAA 2006-6753)}. In: AIAA/AAS
  Astrodynamics Specialist Conference and Exhibit. Guidance, Navigation, and
  Control and Co-located Conferences, American Institute of Aeronautics and
  Astronautics, USA.

\bibitem[{von {Zeipel}(1916, 1917, 1918)}]{vonZeipel1916}
von {Zeipel}, H., July 1916, 1917, 1918. {Research on the motion of minor
  planets (recherches sur le mouvement des petites plan\`etes)}. {NASA
  Translation: NASA TT F-9445. (1965).}
\newline\urlprefix\url{https://archive.org/details/nasa_techdoc_19650019998}

\end{thebibliography}
\end{document}